\documentclass{ifacconf}

\usepackage{natbib}
\usepackage{graphicx}
\usepackage{mathtools}
\usepackage{enumitem}
\usepackage{tikz}
\usepackage{amsmath,amsfonts,amssymb}
\usepackage{rotate}

\DeclareMathOperator{\id}{Id}

\begin{document}

\begin{frontmatter}

\title{YALTAPy and YALTAPy\_Online: Python toolboxes for the $H_\infty$-stability analysis of classical and fractional systems with commensurate delays} 

\author[First]{Hugo Cavalera} 
\author[Second]{Jayvir Raj} 
\author[Third]{Guilherme Mazanti}
\author[Third]{Catherine Bonnet}

\address[First]{Universit\'e Paris-Saclay, CNRS, CentraleSup\'elec, Inria, Laboratoire de signaux et syst\`emes, 91190 Gif-sur-Yvette, France}
\address[Second]{Universit\'e Paris-Saclay, CNRS, CentraleSup\'elec, Inria, Laboratoire de signaux et syst\`emes, 91190 Gif-sur-Yvette, France and IPSA, 63 boulevard de Brandebourg, 94200 Ivry-sur-Seine, France}
\address[Third]{Universit\'e Paris-Saclay, CNRS, CentraleSup\'elec, Inria, Laboratoire de signaux et syst\`emes, 91190 Gif-sur-Yvette, France (e-mail: \texttt{firstname}.\texttt{lastname}@inria.fr)}

\begin{abstract}
The aim of this paper is to give a presentation of the Python toolbox YALTAPy dedicated to the stability study of standard and fractional delay systems as well as its online version YALTAPy\_Online. Both toolboxes are derived from YALTA whose functionalities will be recalled here. Examples will be given to show how these toolboxes may be used.
\end{abstract}


\end{frontmatter}

\section{Introduction}
The stability analysis of linear delay systems is a deep research subject which has been the source of countless studies in time and frequency domains since the 1950's, such as \cite{Pontryagin1955, BC63, Niculescu2004, Richard2003, Walton1987, Niculescu2014}. More recently, the literature involved the stability study of fractional linear systems with delays, as fractional systems successfully model complex phenomena in a compact manner.

Many results have been obtained in both settings concerning stability characterization (asymptotic/exponential stability, $H_\infty$-stability) of retarded delay systems. However, much remains to be done for neutral delay systems. Even though stability characterization (by sufficient or necessary and sufficient conditions) is a central issue, several related questions such as easy-to-check stability conditions, precise location of stable/unstable poles, behavior of stable/unstable poles when the delay varies, are of interest for a full and efficient study of practical situations.  

In parallel of theoretical and numerical studies, e.g.\ \cite{Olgac2004, GM12, hc06, fbon12}, toolboxes have been developed to facilitate and popularize the use of the obtained results both in the academic and industrial environments. We cite here  \emph{Quasi-Polynomial Mapping Based Rootfinder} (QPmR, \cite{QPmR}), \emph{Tool for Robust Analysis and Characteristic Equations of Delay Differential Equations} (TRACE-DDE, \cite{trace}), \emph{Bifurcation analysis of delay differential equations} (DDE-Biftool, \cite{ddebiftool}), \emph{Yet another LTI TDS Algorithm} (YALTA, \cite{afb-tds13}), \emph{Partial pole placement via delay action} (P3$\delta$, \cite{Boussaada2020Partial, Boussaada2021New}), and \emph{Parallel Processing Delay Margin Finder} (parDMF, \cite{Ramirez2021Scalable}). The references \cite{afb-tds13, afbn-springer14} present the functionalities of the Matlab toolbox YALTA and briefly describe its positioning relative to earlier toolboxes. We also refer to \cite{Sipahi_book2019} and \cite{Pekar2018} for a recent overview on effective methods for the stability analysis of delay systems where toolboxes are also presented. 

The aim of this paper is to present YALTAPy, a Python version of the Matlab toolbox YALTA, as well as YALTAPy\_Online, which is an online version of YALTAPy. Note that YALTAPy is an open source software.

In Section~\ref{Mainfeatures} we will describe the main functionalities of YALTAPy (i.e., briefly recall those of YALTA). Section~\ref{sec:yaltaonline} will detail a bit the features of YALTAPy\_Online. Then examples of use of YALTAPy and YALTAPy\_Online will respectively be given in Sections \ref{sec:expl-yaltapy} and \ref{sec:expl-yaltapyonline}. 

\section{Main features of YALTAPy}
\label{Mainfeatures}

YALTAPy and YALTAPy\_Online consider delay systems described by transfer functions of
the type:
\begin{equation}
\label{eq01}
\displaystyle
\textstyle G(s)= \frac{t(s) +  \sum_{\kappa=1}^{N'}t_\kappa(s)e^{-\kappa \tau s} }
{p_0(s)+ \sum_{k=1}^N
p_k(s)e^{-k \tau s}}=\frac{n(s)}{d(s)},
\end{equation}
where $\tau > 0$ is the nominal delay, $\alpha \in (0, 1]$ is the fractional exponent, and, for every $k \in \{0, \dotsc, N\}$, $p_k(s)$ is a polynomial in the variable $s^\alpha$. We let $n$ denote the degree of $p_0$ and we assume that $n \ge \deg t$, $n \ge \deg t_\kappa$ for every $\kappa \in \{1, \dotsc, N'\}$, and $n \geq \deg p_k$ for every $k \in \{1, \dotsc, N\}$, i.e., the system is either of retarded or neutral type. 

We refer to \cite{BC63, fbon10, Partington2004} for an analysis of the position of chains of poles of such systems:
\begin{enumerate}
\item If $\deg p_0 > \deg p_k$ for all $k \in \{1, \dotsc, N\}$, then there are only chains of retarded type; 
\item If $\deg p_0=\deg p_N$, then there are only chains of neutral type;
\item If $\deg p_0=\deg p_k > \deg p_N$, for some $k \in \{1, \dotsc, N-1\}$, then there are chains of both neutral and retarded types.
\end{enumerate}
We recall that a neutral chain is asymptotic to a vertical axis in the complex plane, while a retarded chain contains poles with arbitrarily large negative real part. The coefficient of the highest degree term of $p_0(s) + \sum_{k=1}^N p_k(s)e^{-k \tau s}$ is given by 
\begin{equation}
\label{eq07}
\textstyle \tilde{c}_d(z)= 1 + \sum_{i=1}^N\alpha_i z^i.
\end{equation}

When dealing with neutral systems, in order to avoid the possibility of an infinite number of zero cancellations between the numerator and denominator of $G$, we make the same Hypothesis (H) than in \cite{afb-tds13,afbn-springer14} in terms of  $\deg t$,  $\deg t_k$ and roots of $\tilde{c}_d$ and $\tilde{c}_n$.

We refer to \cite{afb-tds13,afbn-springer14} for a complete description of the functionalities of YALTA and just precise below that YALTAPy and YALTAPy\_Online will give:
\begin{itemize}
\item the type of delay system: retarded or neutral
\item the position of asymptotic axes for a neutral system
\item the stability property when the delay is zero
\item stability windows
\item a root locus showing the displacement of the poles when the delay varies between two chosen values
\item a Pad\'e-2 approximation of the system (only for standard systems).
\end{itemize}

Let us now describe how inputs should be provided to YALTAPy functions. The quasi-polynomial $d(s)$ may be rewritten 
\begin{equation}
\label{eq:delta}
\textstyle d(s) = P_0(s) + \sum_{k=1}^M P_k(s) e^{-n_k \tau s},
\end{equation}
where, for every $k \in \{0, \dotsc, M\}$, $n_k$ is a positive integer and $P_k(s) = p_{n_k}(s)$. In other words, we eliminate from the denominator of \eqref{eq01} all delays $k \tau$ corresponding to values of $k$ for which $p_k$ is identically zero, and $M \in \{1, \dotsc, N\}$ denotes thus the number of polynomials among $p_1, \dotsc, p_N$ that are not identically zero.

The data of the system corresponding to \eqref{eq:delta} is hence described by the coefficients of the polynomials $P_0, \dotsc, P_M$, the exponent $\alpha$, the delay $\tau$, and the integer multiples $n_1, \dotsc, n_M$ of the delay $\tau$ that appear in the system.

In order to input the coefficients of $P_0, \dotsc, P_M$ to YALTAPy, one should provide the $(M+1) \times (n+1)$ matrix
\[
\begin{pmatrix}
p_{0, n} & p_{0, n-1} & \cdots & p_{0, 0} \\
p_{1, n} & p_{1, n-1} & \cdots & p_{1, 0} \\
\vdots & \vdots & \ddots & \vdots \\
p_{M, n} & p_{M, n-1} & \cdots & p_{M, 0} \\
\end{pmatrix},
\]
i.e., the matrix whose $k$th row (starting from $k = 0$) describes the coefficients of $P_k$ in descending order of power. This matrix should be represented in Python as a $2$-dimensional \texttt{numpy} array (see \cite{harris2020array} for details on \texttt{numpy}). Concerning the integer multiples $n_1, \dotsc, n_M$ of the delay $\tau$, they should be input as a vector of size $M$,
\[
\begin{pmatrix}
n_1 & \cdots & n_M \\
\end{pmatrix},
\]
which should be represented as a $1$-dimensional \texttt{numpy} array of integers. Finally, $\tau$ and $\alpha$ should be input to YALTAPy as floating-point numbers or integers. Examples of use of YALTAPy are provided in Section~\ref{sec:expl-yaltapy} below.

\section{YALTAPy\_Online}
\label{sec:yaltaonline}

The team working on YALTA has recently developed an interface for facilitating the use of YALTAPy, based on a Python's \emph{Jupyter Notebook}, an open document format which can contain live code, equations, visualizations, and text. Our Jupyter Notebook implements a friendly graphical user interface for YALTAPy\_Online thanks to interactive widgets from Python's \texttt{ipywidgets} module. The next aim is to facilitate the use of YALTAPy\_Online by making it available with no need of installation thanks to the \emph{Binder} service \citep{project_jupyter-proc-scipy-2018}, which allows the creation of personalized computing environments directly from a \emph{git} repository. The \emph{Binder} service is free to use and is powered by \emph{BinderHub}, an open-source tool that deploys the service in the cloud.

\begin{figure}[ht]
\centering
\includegraphics[width=0.85\columnwidth]{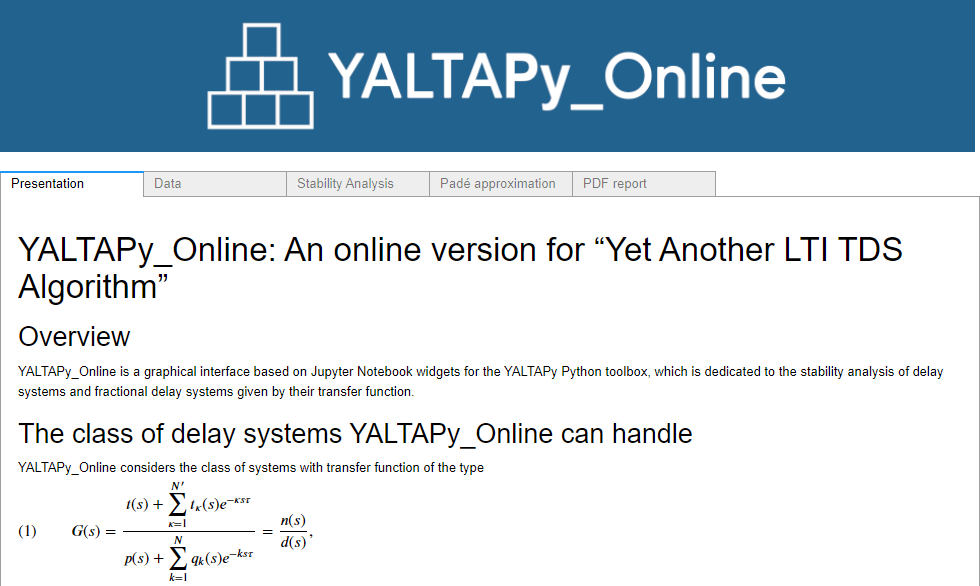}
\caption{Main screen of YALTAPy\_Online.}
\label{fig:online}
\end{figure}

YALTAPy\_Online starts by a presentation screen (see Fig.~\ref{fig:online}) recalling the class of systems that it handles and its main features. It is divided into five tabs: the first one is the presentation screen, the second one is dedicated to the input of the data, the third one provides features of stability analysis, stability windows, and root locus, the third one provides a tool for computing Padé approximations, and the last one allows the user to obtain a PDF report with the obtained results. A detailed example of the use of YALTAPy\_Online is provided in Section~\ref{sec:expl-yaltapyonline}.

\section{Examples of use of YALTAPy}
\label{sec:expl-yaltapy}

In this section, we present two examples illustrating the use of YALTAPy for the stability analysis of time-delay systems. The first one, described in Section~\ref{sec:expl-yaltapy-1}, considers a time-delay system of order $3$ with a single delay, whereas the second one, described in Section~\ref{sec:expl-yaltapy-2}, illustrates how to use YALTAPy to study the stability of a fractional system.

\subsection{A third-order system}
\label{sec:expl-yaltapy-1}

Consider the time-delay system given by
\begin{equation}
\label{eq:sys-expl-1}
\dot x(t) = A_0 x(t) + A_1 x(t - \tau),
\end{equation}
where $x(t) \in \mathbb R^n$, $A_0$ and $A_1$ are $n \times n$ matrices with real coefficients, and $\tau > 0$ is the delay. The stability of such a system can be studied through the roots of its characteristic equation
\[
d(s) = \det(s \id - A_0 - A_1 e^{-s \tau}),
\]
where $\id$ denotes the $n \times n$ identity matrix.

In this example, we consider $n = 3$, $\tau = 3$, and
\[
\textstyle A_0 = \begin{pmatrix}
0 & 1 & -1 \\
1 & -\tfrac{1}{2} & -\tfrac{1}{4} \\
-1 & 0 & -\tfrac{4}{5}
\end{pmatrix}, \qquad A_1 = \begin{pmatrix}
-3 & -2 & \tfrac{1}{2} \\
-1 & 0 & \tfrac{1}{5} \\
\tfrac{2}{3} & \tfrac{1}{3} & 0
\end{pmatrix}.
\]
In this case, a straightforward computation shows that
\begin{equation}
\label{eq:delta-1}
\begin{split}
d(s) = {} & \textstyle s^3 + \frac{13}{10} s^2 - \frac{8}{5} s - \frac{31}{20} + \textstyle {} + \left(3 s^2 + \frac{163}{20} s + \frac{323}{60}\right) e^{-\tau s}   \\
& \textstyle {} - \left(\frac{12}{5}s + \frac{173}{60}\right) e^{-2 \tau s} + \frac{7}{30} e^{-3 \tau s}.
\end{split}
\end{equation}

In order to input the characteristic equation \eqref{eq:delta-1} into YALTAPy, we define the \texttt{numpy} array
\begin{verbatim}
P = np.array([[1, 13/10,  -8/5 ,  -31/20],
              [0,     3, 163/20,  323/60],
              [0,     0,  -12/5, -173/60],
              [0,     0,      0,    7/30]])
\end{verbatim}
(assuming that \texttt{numpy} was previously imported using \texttt{import numpy as np}). Similarly, we input the vector \texttt{n} of the multiples of the delays and the values of $\tau$ and $\alpha$ by
\begin{verbatim}
n = np.array([1, 2, 3]), tau = 3, alpha = 1
\end{verbatim}

\subsubsection{Stability analysis}

With the previous definitions, one studies the stability of \eqref{eq:sys-expl-1} by using YALTAPy's \texttt{thread\_\allowbreak{}analysis} function. After importing YALTAPy using \texttt{import yaltapy as yp} and defining \texttt{P}, \texttt{n}, \texttt{tau}, and \texttt{alpha} as above, the \texttt{thread\_\allowbreak{}analysis} function is called by
\begin{verbatim}
res = yp.thread_analysis(P, n, alpha, tau)
\end{verbatim}

After its execution, \texttt{thread\_\allowbreak{}analysis} returns a dictionary with information on the stability of the system. In this example, the variable \texttt{res} contains the following data:
\begin{verbatim}
Type: Retarded
AsympStability: There is (are) 4 unstable
pole(s) in the right half-plane
RootsNoDelay: [-3.07585975+0.j
         -0.61207012+0.10043131j
         -0.61207012-0.10043131j]
RootsChain: Roots chains only computed for
neutral systems
CrossingTable:
 [[0.38196552 1.86170973 3.37495433 1.        ]]
ImaginaryRoots:
[[ 0.38196552  3.37495433  1.        ]
 [ 0.38196552 -3.37495433  1.        ]
 [ 2.24367525  3.37495433  1.        ]
 [ 2.24367525 -3.37495433  1.        ]]
\end{verbatim}

Hence, YALTAPy has identified that the given system is of retarded type and that it has 4 poles in the right half-plane, yielding the conclusion that, with the given parameters, the system is unstable. YALTAPy has also computed the characteristic roots of the system without delay, i.e., when $\tau = 0$. In this case, the characteristic equation of the system reduces to the polynomial $d(s) = s^{3} + \frac{43}{10} s^{2} + \frac{83}{20}s + \frac{71}{60}$, whose three complex roots are those provided in YALTAPy output \texttt{RootsNoDelay}.

Since the system is of retarded type, it presents no chains of roots asymptotic to a vertical line, a fact that is recalled in the output \texttt{RootsChain}. The outputs \texttt{CrossingTable} and \texttt{ImaginaryRoots} describe the behavior of the system as the delay increases from $0$ to its nominal value. For instance, \texttt{ImaginaryRoots} states that, at $\tau \approx 0.3820$, two roots cross the imaginary axis, at the values $\pm 3.3750 i$, and the direction of crossing is from the left to the right (value $1$ in the third column), and, at $\tau \approx 2.2437$, two other roots cross the imaginary axis, once again at the values $\pm 3.3750 i$ and from the left to the right. As for the \texttt{CrossingTable}, its first column gives the first value of the delay for which a specific zero crosses the axis, its third column gives the frequency $\omega$ of crossing of a zero, its second column is $\frac{2\pi}{\omega}$, and its fourth column provides the crossing direction.

\subsubsection{Stability windows}

YALTAPy provides a function to obtain stability windows for the system under consideration, called \texttt{thread\_\allowbreak{}stability\_\allowbreak{}windows}. After defining \texttt{P}, \texttt{n}, \texttt{tau}, and \texttt{alpha} as above, this function is called by
\begin{verbatim}
res = yp.thread_stability_windows(P, n,
      alpha, tau)
\end{verbatim}

As for \texttt{thread\_\allowbreak{}analysis}, the function \texttt{thread\_\allowbreak{}stability\_\allowbreak{}windows} also returns a dictionary with information on stability windows for the system. In the present example, in addition to outputs already given by \texttt{thread\_\allowbreak{}analysis}, the variable \texttt{res} contains the following data:
\begin{verbatim}
StabilityWindows:
 [[0.         0.38196552 3.        ]
  [1.         0.         0.        ]]
NbUnstablePoles:
 [[0.         0.38196552 2.24367525 3.        ]
  [0.         2.         4.         4.        ]]
\end{verbatim}

The outputs \texttt{StabilityWindows} and \texttt{NbUnstablePoles}, in this example, are interpreted as follows: at $\tau = 0$, the system is stable (value $1$ at \texttt{StabilityWindows}) and has no poles in the right half-plane (value $0$ at \texttt{NbUnstablePoles}). At $\tau \approx 0.3820$, the system becomes unstable (value $0$ at \texttt{StabilityWindows}) and has $2$ poles in the right half-plane (value $2$ at \texttt{NbUnstablePoles}). At $\tau \approx 2.2437$, two additional poles appear in the complex right half-plane, and the system remains unstable and with $4$ poles in the right half-plane until the nominal delay value $\tau = 3$.

\begin{figure}[ht]
\centering
\resizebox{0.85\columnwidth}{!}{\input{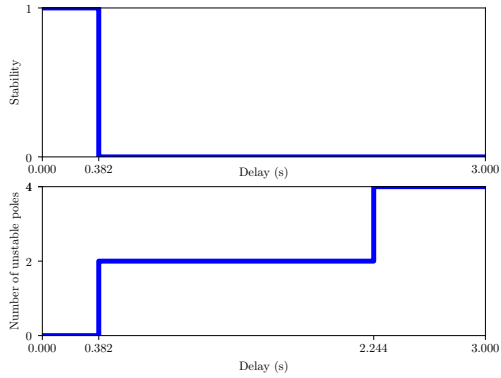}}
\caption{Stability windows for the example from Section~\ref{sec:expl-yaltapy-1}.}
\label{fig:expl1-stab-win}
\end{figure}

In addition to those outputs, \texttt{thread\_\allowbreak{}stability\_\allowbreak{}windows} also plots the data in \texttt{Stability\allowbreak{}Windows} and \texttt{Nb\allowbreak{}Unstable\allowbreak{}Poles} by setting its optional input argument \texttt{plot} to \texttt{True} (default). The plot for this example is provided in Fig.~\ref{fig:expl1-stab-win}.

\subsubsection{Root locus}

The root locus of the roots of a system that eventually become unstable is computed using YALTAPy's function \texttt{thread\_\allowbreak{}root\_\allowbreak{}locus}. For systems of neutral type, this function only runs if the system with nominal delay possesses chains asymptotic to vertical axes located in the open left half-plane.Similarly to the previous functions, this function is called by
\begin{verbatim}
res = yp.thread_root_locus(P, n, alpha, tau)
\end{verbatim}
and it returns a dictionary with information on the root locus. In the present example, the dictionary \texttt{res} returned by the function contains, in addition to outputs also provided by the previous functions, the outputs
\begin{verbatim}
UnstablePoles:
[0.08427953+2.50094915e+00j
 0.08427953-2.50094915e+00j
 1.01660266+2.18573355e-13j]
PolesError: [1.e-11+0.j 1.e-12+0.j 1.e-10+0.j]
\end{verbatim}
The dictionary also contains an entry \texttt{RootLocus} containing information on the movement of the roots, but since it is too long, it is not shown here.

\begin{figure}[ht]
\centering
\resizebox{0.8\columnwidth}{!}{\input{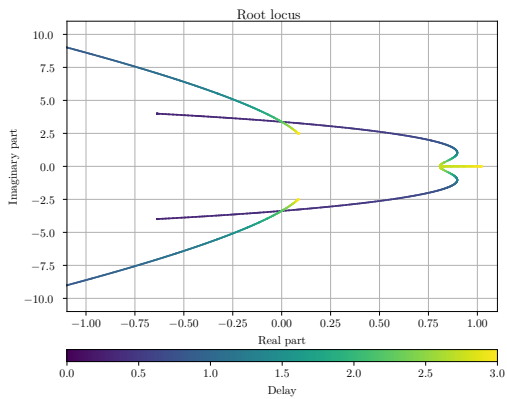}}
\caption{Root locus for the example from Section~\ref{sec:expl-yaltapy-1}.}
\label{fig:expl1-root-locus}
\end{figure}

The output \texttt{UnstablePoles} provides the location of the $4$ unstable poles at the nominal value of the delay $\tau = 3$ (we have a real pole of multiplicity $2$ in this example) and the respective estimated numerical errors in \texttt{PolesError}. If the optional input argument \texttt{plot} is set to \texttt{True}, the root locus is graphically displayed, as shown in Fig.~\ref{fig:expl1-root-locus}, where we observe the behavior of the four roots of the system that become unstable as the delay increases from $0$ to $3$.

\subsection{A fractional system}
\label{sec:expl-yaltapy-2}

As a second example, we consider the fractional system described by its characteristic equation
\[
\textstyle d(s) = s^{\frac{1}{2}} + \frac{15}{2} - \left(\frac{1}{6} s^{\frac{1}{2}} + \frac{20}{7}\right) e^{-\tau s} + \frac{80}{11} e^{-3 \tau s},
\]
with $\tau = 2.5$. Such a system is input into YALTAPy by defining the variables \texttt{P}, \texttt{n}, \texttt{tau}, and \texttt{alpha} as
\begin{verbatim}
P = np.array([[   1,  15/2],
              [-1/6, -20/7],
              [   0,  80/11]])
n = np.array([1, 3]), alpha = 0.5, tau = 2.5
\end{verbatim}

\subsubsection{Stability analysis}

Similarly to the example in Section~\ref{sec:expl-yaltapy-1}, we run \texttt{thread\_\allowbreak{}analysis} in the present example, obtaining as a result
\begin{verbatim}
Type: Neutral
AsympStability: There is (are) 4 unstable
pole(s) in the right half-plane
RootsNoDelay: []
RootsChain: [-0.71670379]
CrossingTable:
[[0.42100844 2.96762299 2.11724512 1.        ]
 [1.96353668 2.32594339 2.70134919 1.        ]]
ImaginaryRoots:
[[ 0.42100844  2.11724512  1.        ]
 [ 0.42100844 -2.11724512  1.        ]
 [ 1.96353668  2.70134919  1.        ]
 [ 1.96353668 -2.70134919  1.        ]]
\end{verbatim}
YALTAPy identifies that the system has a single chain of poles asymptotic to the vertical line with real part $-0.7167$, indicates that the system is unstable for the nominal value of the delay, and identifies the crossings of the imaginary axis as the delay increases from $0$ to $2.5$, providing the value of the delay at which the crossing occurs and the position and direction of crossing.

\subsubsection{Stability windows}

Using YALTAPy to compute stability windows, we obtain the following output (only fields not already returned by \texttt{thread\_\allowbreak{}analysis} are shown)
\begin{verbatim}
StabilityWindows:
[[0.         0.42100844 2.5       ]
 [1.         0.         0.        ]]
NbUnstablePoles:
[[0.         0.42100844 1.96353668 2.5       ]
 [0.         2.         4.         4.        ]]
\end{verbatim}

\subsubsection{Root locus}

\begin{figure}[ht]
\centering
\resizebox{0.8\columnwidth}{!}{\input{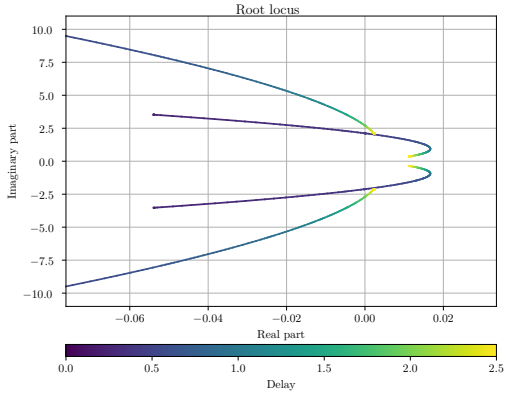}}
\caption{Root locus for the example from Section~\ref{sec:expl-yaltapy-2}.}
\label{fig:expl2-root-locus}
\end{figure}

Finally, we use YALTAPy to compute the root locus of the present example, obtaining the following output (outputs already present in the previous functions are omitted, as well as the output \texttt{RootLocus}, too long to be shown here)
\begin{verbatim}
UnstablePoles:
[0.00224791+2.12366734j
 0.00224791-2.12366734j
 0.01138514-0.36488849j
 0.01138514+0.36488849j]
PolesError: [1.e-10+0.j 1.e-13+0.j
             1.e-12+0.j 1.e-11+0.j]
\end{verbatim}
The corresponding root locus is presented in Fig.~\ref{fig:expl2-root-locus}, where we observe in particular the crossings already known from the previous computation of stability windows.

\subsection{Padé approximation}

YALTAPy can compute Padé approximation of functions of the form $(p_0(s) + \sum_{k=1}^N p_k(s) e^{-k s \tau})/(s + 1)^\delta$ where $p_0(s), \dotsc, p_k(s)$ are polynomials in the variable $s$ following the same assumptions as before and $\delta$ is an integer greater than the degree of $p_0$. We illustrate this functionality of YALTAPy in the case $\delta = 2$ and
\begin{equation}
\label{eq:Pade}
d(s) = s + 2 - e^{-s}.
\end{equation}
To input such a system into YALTAPy, we define
\begin{verbatim}
P = np.array([[1, 2], [0, -1]]), 
n = np.array([1]), tau = 1, delta = 2
\end{verbatim}
Padé approximations can be computed through YALTAPy's \texttt{compute\_pade} function, which, in addition to the above variables, takes as arguments also the mode of approximation (``ORDER'' or ``NORM'') and a parameter \texttt{mod\_arg} whose interpretation depends on the mode: \texttt{mod\_arg} is the order of approximation in mode ``ORDER'' and the maximum error in $H_\infty$ norm for the approximation in mode ``NORM''. For instance, a first-order Padé approximation of \eqref{eq:Pade} is obtained by
\begin{verbatim}
res = yp.compute_pade(P, delta, tau, n,
         mod_arg = 1, mode = "ORDER")
\end{verbatim}
The variable \texttt{res} will contain the data of the approximation: \texttt{res.num\_approx} and \texttt{res.den\_approx} contain the coefficients of the numerator and denominator, respectively, of the Padé approximation, while \texttt{res.error\_norm} and \texttt{res.pade\_order} contain, respectively, the $H_\infty$ norm of the approximation error and the order of approximation. After running the present example, these four variables contain, respectively, the values
\begin{verbatim}
array([ 1.,  9., 45., 79., 54., 12.])
array([ 1., 10., 42., 88., 97., 54., 12.])
0.032917024033214635
1
\end{verbatim}

\section{Example of use of YALTAPy\_Online}
\label{sec:expl-yaltapyonline}

\begin{figure}[ht]
\centering
\includegraphics[width=0.85\columnwidth]{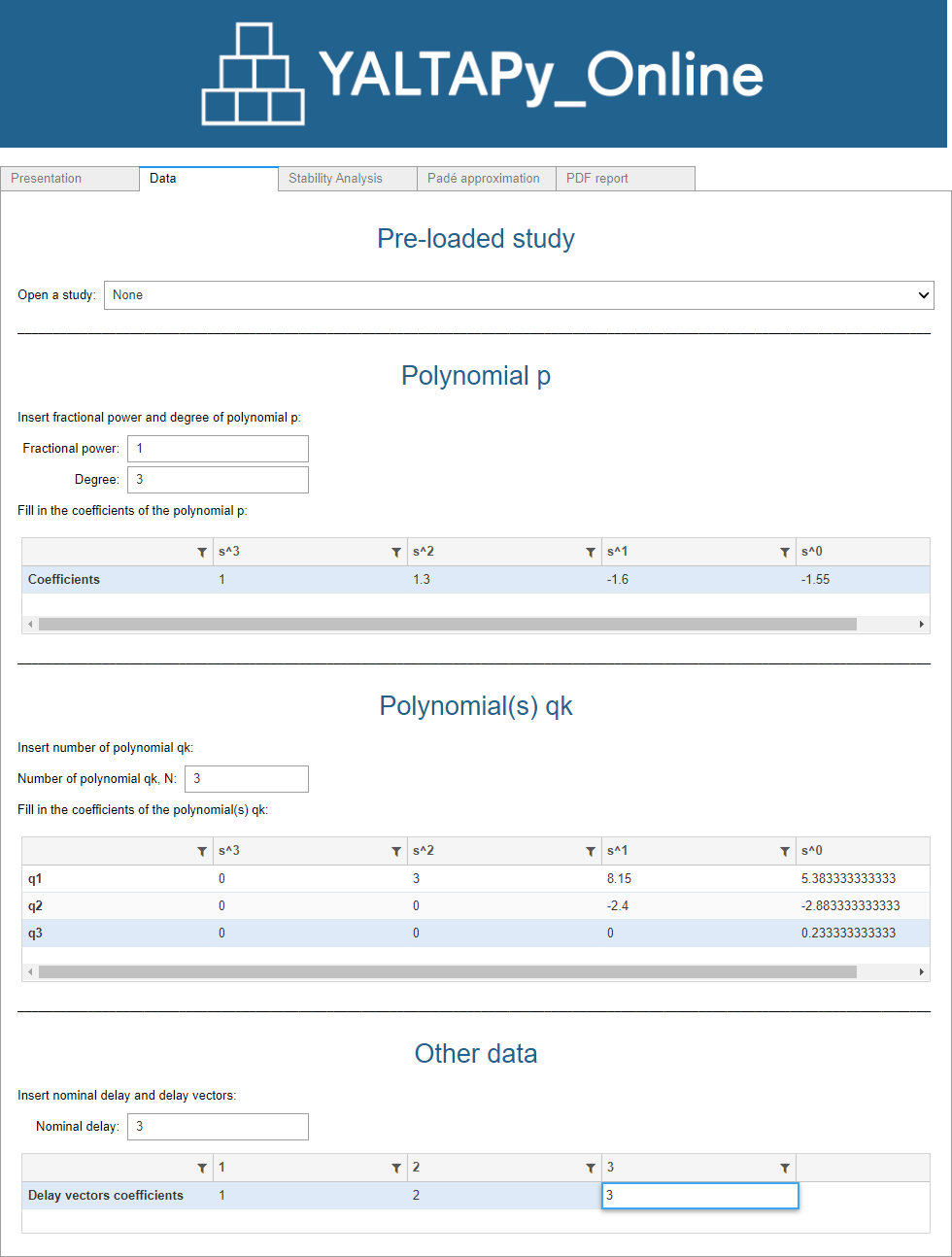}
\caption{Data input in YALTAPy\_Online.}
\label{fig:online-input}
\end{figure}

Let us now illustrate the use of YALTAPy\_Online through the example from Section~\ref{sec:expl-yaltapy-1} above. The first step is to input the data of the problem in the tab ``Data'' (see Fig.~\ref{fig:online-input}). The polynomials without delays and with delays are entered separately. After selecting the fractional power of $s^\alpha$ (which is $\alpha = 1$ in this example) and degree of the polynomial $P_0$ without delay terms ($3$ in this example), the user enters its coefficients in the corresponding field. The user is then prompted for the number of delayed terms and for the coefficients of the corresponding polynomials. Finally, the last part is dedicated to the nominal delay and the vector containing the integers multiplying the delay.

\begin{figure}[ht]
\centering
\includegraphics[width=0.85\columnwidth]{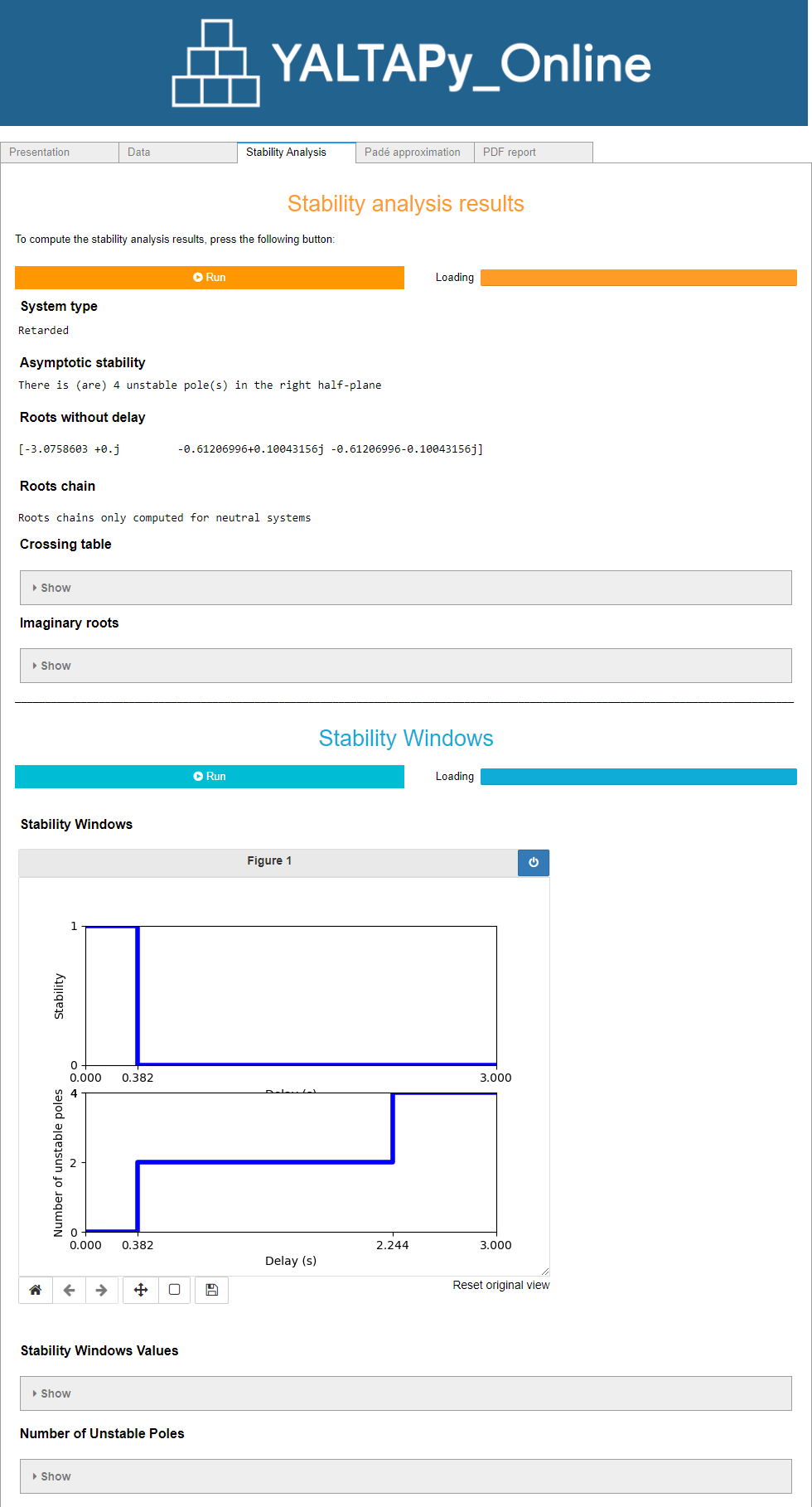}
\caption{Stability analysis in YALTAPy\_Online (root locus omitted for simplicity).}
\label{fig:online-analysis}
\end{figure}

After entering the data, the user goes to the tab ``Stability Analysis'' to obtain YALTAPy\_Online's study of the system. As YALTAPy, YALTAPy\_Online will perform a stability analysis (corresponding to YALTAPy's function \texttt{thread\_\allowbreak{}analysis}), compute stability windows (corresponding to \texttt{thread\_\allowbreak{}stability\_\allowbreak{}windows}), and compute root locus of unstable roots (corresponding to \texttt{thread\_\allowbreak{}root\_\allowbreak{}locus}). Those analyses are started by the respective ``Run'' buttons on the page and all results are displayed on the screen. Long outputs are automatically hidden, but can be shown through the corresponding ``Show'' buttons. Graphs are interactive, and the user can zoom in parts of the graph or export it as an image. The screen obtained after running all available analyses for the example from Section~\ref{sec:expl-yaltapy-1} is provided in Fig.~\ref{fig:online-analysis} (the root locus was omitted due to space constraints).

\end{document}